\documentstyle{article}
\newcommand{\be}{\begin{equation}}
\newcommand{\ee}{\end{equation}}
\newcommand{\ap}{a^{\dagger}}
\newcommand{\tlambda}{\tilde{\lambda}}
\newcommand{\tla}{\tilde{\lambda}}

\newcommand{\tnu}{\tilde{\nu}}

\newcommand{\s}{a^{\dag}}

\newcommand{\cc}{|c,\nu_0 \, \rangle }
\newcommand{\nn}{| c ,\nu_0 + n ) }
\newcommand{\no}{| c ,\nu_0 + n \, \rangle }
\newcommand{\la}{\lambda}
\title{Representation Theory of Generalized Deformed Oscillator
Algebras \thanks{Presented at the 5th Colloquium ''Quantum Groups and
Integrable Systems'', Prague, 20-22 June 1996.}}
\author{Christiane Quesne \thanks{Directeur de recherches FNRS; 
e-mail: cquesne@ulb.ac.be}  \hspace{1mm} and \hspace{1mm} 
{Nicolas Vansteenkiste} \thanks{e-mail: nvsteen@ulb.ac.be}}
\date{}
\begin{document}
\maketitle
\begin{abstract}
The representation theory of the generalized deformed                          
oscillator algebras (GDOA's) is developed. GDOA's are generated by 
the four operators $ \{ {\bf 1}, a, a^{\dag} , N \} $. Their 
commutators and Hermiticity properties are those of the boson 
oscillator algebra, except for $ [ a, a^{\dag} ]_q = G(N)$, 
where $[a,b]_q = a b - q \, b a $ and $ G(N) $ is a Hermitian, 
analytic function.
The unitary irreductible representations are obtained 
by means of a Casimir operator $C$ and the semi-positive operator 
$ a^{\dag} a$. They may belong to one out of four classes: 
bounded from below (BFB), bounded from above (BFA), 
finite-dimentional (FD), unbounded (UB).
Some examples of these different types of unirreps 
are given.
\end{abstract}
\section{Introduction}
Since the pioneering works of Arik and Coon~\cite{arik},
Kuryshkin, Biedenharn, and Macfarlane~\cite{kuryshkin}, 
many forms of deformed oscillator algebras 
(DOA's) have been considered. They have played an important role in the 
construction of $q$-deformed Lie algebras 
and have found various applications to physical problems. 

The necessity to introduce some order in the rich and varied choice of 
deformed commutation relations did however appear soon and various 
classification schemes were therefore proposed
~\cite{jannussis,daska91,katriel}.
The representations of some forms of DOA's have 
been investigated~\cite{kulish}, which brought out 
the existence of additional non-Fock-space representations. 

In this communication, the general representation theory of 
generalized DOA's will be considered \cite{gdoarep}.
The representation theory of generalized deformed oscillator algebras 
is developed in Sec.~2, and illustrated on some examples in Sec.~3. 
Finally, Sec.~4 contains the conclusion.
\section{Representation theory of generalized deformed oscillator
algebras}
1. A  generalized deformed oscillator algebra (GDOA) is 
generated by the operators
$ \{ {\bf 1}, a, \s , N \} $ satisfying the Hermiticity conditions
$ ( \s )^{\dag} = a $, $ N^{\dag} = N $, and the commutation
relations
\be
[N,a] = -a \; \; \; \; \; \; \; \; \; \; \; \;  
            \; \; \; \; \; \; \; \; \; \; \; \; [N , \s ] = \s \: , 
\ee
\be
[a, \s ]_q = G(N) \: \label{2} , 
\ee
where $ [a,b]_q = a b - q \, ba $, $ q \in R$, is a quommutator, and
the deformation function $ G(N) $ is a Hermitian, analytic function.\\
The structure function, $ F(N) $, is defined as the
solution of the functional equation
\be
F(x+1) - q \, F(x) = G(x) \: , \; \;  x \in R \; \; 
                                {\rm and} \; \; F(0) = 0 \: . 
\ee
 
2. Contrary to some definitions of GDOA's where
both a commutation and an anticommutation relations are
assumed~\cite{daska91}, we only keep a commutator. As explained in
Refs.~\cite{oh}, this choice leads to the existence of a
Casimir operator $C$ that can be constructed by means of the structure
function:
\be
C = q^{-N} (F(N) - \s a ) = q^{-(N+1)} (F(N+1) - a \s) \, . 
\ee
This implies an algebraic relation between $ \s a $ ($ a\s $) and
$N$, $C$
\be
\s a = F(N) - q^N C \; , \; \; a\s = F(N+1) - q^{N+1} C \,. \label{5}
\ee
Another useful Casimir operator is $ U= {\rm e}^{{\rm i} 2 \pi N} $.

3. The spectrum of $N$ is discrete. This can be proved by the same 
technique as that used for parabosons in \cite{jordan}. The Casimir 
operator $ U $ is 
unitary, so that in a given unirrep, it possesses a fixed eigenvalue 
of the form $ {\rm e}^{{\rm i} 2 \pi \nu_0} $, $\nu_0 \in R$. On 
the other hand, the eigenvalues of $ U $ can be determined from those 
of the Hermitian operator $N$. The spectral mapping theorem leads to 
eigenvalues for $U$ of the form 
$ {\rm e}^{{\rm i} 2 \pi x} $, where $ x \in R$ are the eigenvalues 
of $N$. The equivalence of the two expressions for the eigenvalues of 
$U$ implies that $ x= \nu_0 + n $, $n \in Z $, in a given unirrep.

The first assumption we make is that the spectrum of $N$ is 
nondegenerate. Then, we suppose the existence of a normalized 
simultaneous eigenvector, $ \cc $, of $C$ and $N$:
\begin{eqnarray}
                   C \cc & = & c \cc \, ,  \\
                   N \cc & = & \nu_0 \cc \, ,  \\
\langle \, c , \nu_0 \cc & = &  1 \, .  
\end{eqnarray}    

4. By the repeated action of $a$ and $\s$ on this vector, it is possible
to construct new eigenvectors of $C$ and $N$. With
\begin{eqnarray}
\nn & \equiv   & ( \s )^n \, \cc \; \; \; \; \; \; \; \, \; {\rm if} 
                                \; \; \; \; \; n > 0   \\
    & \equiv   & ( a )^{-n} \, \cc \; \; \; \; \; \; \; \; {\rm if}
                                \; \; \; \; \; n<0 \; , 
\end{eqnarray}
we get
\begin{eqnarray}
C \, \nn    & = & c \, \nn \, ,  \\
N \, \nn    & = & ( \nu_0 + n ) \, \nn \, ,  
\end{eqnarray}
if $ \nn $ is non-vanishing.

These vectors are also eigenvectors of the semi-positive operators
$ \s a $ and $a \s $:
\begin{eqnarray}
\s a \, \nn & = & \la_n \, \nn \, , \\
a \s \, \nn & = & \mu_n \, \nn \, ,  
\end{eqnarray}
with, according to (\ref{5}),
\begin{eqnarray}
\la_n & = & F(\nu_0+n) - q^{\nu_0+n} c \, , \\
\mu_n & = & \la_{n+1} \, .
\end{eqnarray}
The existence of $ \cc $ and the positiveness of $ \s a $ impose a
condition on $ \la_0 $:
\be
\la_0 = F(\nu_0) - q^{\nu_0} c \geq 0 \, \label{la} .
\ee

As long as they exist, the vectors $ \nn $ can be normalized:
\begin{eqnarray}
\no & \equiv   & \left( \prod_{i=1}^{n} \la_i \right)^{-1/2} \, 
             ( \s )^n  \, \cc   \; \; \; \; \; \; \;   {\rm if} 
                           \; n = 0 , 1 , 2 , \ldots  \label{rpo} \\
    & \equiv   & \left( \prod_{i=0}^{|n|-1} \la_{-i} \right)^{-1/2} \,
                        ( a )^{-n} \, \cc    \; {\rm if}
                            \; n=0, -1 , -2 , \ldots \, \label{rne} .
\end{eqnarray}
In the representation spanned by these vectors, the matrix elements 
of the GDOA's generators are
\begin{eqnarray}
\langle \, c, \nu_0 + m | \, C \, | c, \nu_0 + n \, \rangle & = & 
                                       c \: \delta_{m,n} \, ,   \\ 
\langle \, c, \nu_0 + m | \, N \, | c, \nu_0 + n \, \rangle & = & 
                               (\nu_0 +n) \: \delta_{m,n} \, ,   \\
\langle \, c, \nu_0 + m | \, a \, | c, \nu_0 + n \, \rangle & = & 
                         \sqrt{\la_n}  \: \delta_{m,n-1}  \, ,   \\
\langle \, c, \nu_0 + m | \, \s \, | c, \nu_0 + n \, \rangle & = & 
                       \sqrt{\la_{n+1}} \: \delta_{m,n+1} \, . 
\end{eqnarray}
The Fock-space representations are characterized by $ c=\nu_0 = 0$ .

5. The existence condition of $ \nn $ (and thus of $ \no $) derives 
from the unitary condition for the representation. The latter 
imposes the positivity of the $ \s a $ eigenvalues, 
as expressed in the following proposition:

{\bf Proposition :}
If there exists some $m_1 \in \{-1, -2, -3, \ldots\}$ such that 
$\lambda_{m_1} < 0$, and $\lambda_n \ge 0$ for $n = 0$, $-1$, 
$\ldots$,~$m_1+1$, then an irreducible representation of a deformed 
oscillator algebra can be unitary only if $\lambda_{n_1} = 0$ for some 
$n_1 \in \{0, -1, \ldots, m_1+1\}$. If there exists some 
$m_2 \in \{2, 3, 4, \ldots \}$ such that $\lambda_{m_2} < 0$, and 
$\lambda_n \ge 0$ for $n = 0$, $1$, $\ldots$,~$m_2-1$, then it can be 
unitary only if $\lambda_{n_2} = 0$ for some 
$n_2 \in \{1, 2, \ldots, m_2-1\}$.

{\it Proof.} In the first part of the proposition, we must have $|c,
\nu_0+m_1) = 0$ as otherwise $\ap a$ would have a negative eigenvalue. This
implies that $a |c, \nu_0+m_1+1) = 0$, which can be achieved in two ways, 
either $|c, \nu_0+m_1+1) = 0$, or $|c, \nu_0+m_1+1) \ne 0$ and 
$\lambda_{m_1+1} = 0$. In the former case, we can proceed in the same way 
and find that at least one of the conditions $\lambda_{m_1+2} = 0$, 
$\lambda_{m_1+3} = 0$, $\ldots$, $\lambda_{-1} = 0$, or 
$|c, \nu_0-1) = 0$ must be satisfied. But the last one is
equivalent to $\lambda_0 = 0$, since $|c, \nu_0) \ne 0$ by hypothesis. This    
concludes the proof of the first part of the proposition. The second part 
can be demonstrated in a similar way by using $a \ap$, and 
$\mu_n = \lambda_{n+1}$ instead of $\ap a$, and~$\lambda_n$.

6. This leads to four classes of unirreps:
\begin{itemize}
\item If there exists some $ n_1 \in \{ 0,-1, -2, \ldots \} $ such that
$ \la_{n_1} = 0 $  and $ \la_n > 0 $ for all 
$  n \in \{n_1+1, n_1+2, \ldots \} $, then $ | c ,\nu_0 + n_1 \, \rangle $
is annihilated by $ a $, and one gets a {\it bounded from below} (BFB) 
unirrep. By repeating the unirrep construction beginning with
$ | c ,\nu_0 + n_1 \, \rangle \equiv  | c ,\tnu_0 \, \rangle $ and 
replacing $ \nu_0 $ by $ \tnu_0 = \nu_0 + n_1 $, $ \la_n $ by
$\tla_n = \la_{n+n_1} $ ($\tla_0 = 0 $), the vector space of the BFB
unirrep is spanned by the vectors defined in
(\ref{rpo}). The eigenvalue of the Casimir operator $C$ is 
$ c= q^{-\tnu_0} F(\tnu_0)$.  
\item If there exists some $ n_2 \in \{ 1, 2, 3 \ldots \} $ such that
$ \la_{n_2} = 0 $  and $ \la_n > 0 $ for all 
$  n \in \{n_2-1, n_2-2, \ldots \} $, then $ | c ,\nu_0 + n_2-1 \, \rangle $
is destroyed by $ \s $, and one gets a {\it bounded from above} (BFA) 
unirrep. By repeating the unirrep construction beginning with
$ | c ,\nu_0 + n_2-1 \, \rangle \equiv  | c ,\tnu_0 \, \rangle $ and 
replacing $ \nu_0 $ by $ \tnu_0 = \nu_0 + n_2 -1 $, $ \la_n $ by
$\tla_n = \la_{n+n_2 -1} $ ($\tla_1 = 0 $), the vector space of the 
BFA unirrep is spanned by the vectors defined in
(\ref{rne}), with $ c= q^{-(\tnu_0+1)} F(\tnu_0+1)$.
\item If there exists some $ n_1 \in \{ 0,-1, -2, \ldots \} $ 
and $ n_2 \in \{ 1, 2, 3 \ldots \} $ such that
$ \la_{n_1} = \la_{n_2} = 0 $ and $ \la_n > 0 $ for all 
$  n \in \{n_1+1, n_1+2, \ldots ,n_2 -1 \} $, then 
$ | c ,\nu_0 + n_1 \, \rangle $ ($ | c ,\nu_0 + n_2-1 \, \rangle $)
is annihilated by $ a $ ($ \s $), and one gets a {\it finite-dimensional}
(FD) unirrep with dimension $ d = p+1 = n_2 - n_1 $. By repeating the 
unirrep construction beginning with 
$ | c ,\nu_0 + n_1 \, \rangle \equiv  | c ,\tnu_0 \, \rangle $ and 
changing $ \nu_0 $ into $ \tnu_0 = \nu_0 + n_1 $, $ \la_n $ into
$\tla_n = \la_{n+n_1} $ ($\tla_0 = \tla_{p+1} = 0 $), the vector space
of the FD unirrep is spanned by the vectors defined in (\ref{rpo}) with 
$ n \leq p $, . The eigenvalue of $C$ is  
$ c= q^{-\tnu_0} F(\tnu_0) = q^{-(\tnu_0+p+1)} F(\tnu_0+p+1)$. \\
The same results are obviously obtained by repeating the unirrep
construction beginning with 
$ | c ,\nu_0 + n_2-1 \, \rangle \equiv  | c ,{\bar \nu}_0 \, \rangle $, 
which leads to a unirrep with vectors defined by (\ref{rne}) and
$ n \geq -p $.
\item If $ \la_n > 0 $ for all $  n \in Z $, one gets an {\it unbounded} 
(UB) unirrep on a vector space spanned by the vectors defined in
(\ref{rpo}) and (\ref{rne}) and $ \nu_0 \in [0,1[ $ (the representations
corresponding to values of $ \nu_0 $ differing by integer steps are
equivalent). 
\end{itemize}
\section{Examples}
To illustrate the general representation theory of GDOA's, we 
consider the
unirreps of three popular GDOA's derived from the boson algebra by the 
recursive minimal-deformation procedure \cite{katriel}.

Note that we consider the case $q < 0$. This case
is omitted in most studies, because the corresponding algebras are 
considered as
deformations of the fermion oscillator algebra. For $ q $ negative,
the parameter $ \nu_0 $ is restricted to integer values so that
$ q^{\nu_0} $ is well defined. 
\subsection{The Arik-Coon oscillator algebra}
It is characterized by $ G(N) = 1 $ \cite{arik}, so that
\be 
F(N)= \frac{q^{N}-1}{q-1} \equiv [N]_q \, ,
\ee
which gives
\be
C = \frac{1-q^{-N}}{q-1} - q^{-N} \s a \, , 
\ee
and
\be
\la_n = \left( \frac{1}{q-1} - c \right) q^{\nu_0+n} - \frac{1}{q-1} \, .
\ee
The different classes of unirreps of this GDOA are listed in  
table 1. The UB unirreps, which diverge for $ q \rightarrow 1^- $, 
are referred to as classically singular representations \cite{aizawa}. 
%
%
\footnotesize
\begin{table}[htb]
\caption{Unirrep classification for the Arik-Coon oscillator algebra. The cases
where $q=1$ and $q=-1$ correspond to the boson and fermion oscillator algebras,
respectively.}
\vspace{.5cm}
\begin{tabular}{lll}          
  \hline\\[-0.2cm] 
  $q$ & Type & Characterization \rule[-0.3cm]{0cm}{0.6cm}\\[0.2cm]
  \hline\\[-0.2cm] 
  $q>1$       & BFB & $\tnu_0\in R$, $c=q^{-\tnu_0} [\tnu_0]_q$, $\tlambda_n = [n]_q$
                       \\[0.2cm]
  $q=1$      & BFB & $\tnu_0\in R$, $c=\tnu_0$, $\tlambda_n = n$ \\[0.2cm]
  $0<q<1$   & BFB & $\tnu_0\in R$, $c=q^{-\tnu_0} [\tnu_0]_q$, $\tlambda_n =
                       [n]_q$\\[0.2cm]
                 & UB   & $0\le\tnu_0<1$, $c\le (q-1)^{-1}$, $\tlambda_n = [\tnu_0+n]_q -
                       c q^{\tnu_0+n}$\\[0.2cm]
  $-1<q<0$ & BFB & $\tnu_0\in Z$, $c=q^{-\tnu_0} [\tnu_0]_q$, $\tlambda_n =
                       [n]_q$\\[0.2cm]
                 & UB  & $0\le\tnu_0<1$, $c=(q-1)^{-1}$, $\tlambda_n=(1-q)^{-1}$ 
                       \\[0.2cm]
  $q=-1$    & FD  & $\tnu_0\in 2 Z$, $p=1$, $c=0$, $\tlambda_n = \left(1 - 
                       (-1)^n\right)/2$ \\[0.2cm]
                & FD   & $\tnu_0\in 2 Z + 1$, $p=1$, $c=-1$, $\tlambda_n = \left(1 - 
                       (-1)^n\right)/2$ \\[0.2cm]
                & UB   & $\tnu_0=0$, $-1<c<-1/2$ or $-1/2<c<0$, \\[0.2cm] 
                &        & $\tlambda_n = (-1)^{n+1}c + \left(1 - (-1)^n\right)/2$ \\[0.2cm]
                & UB   & $0\le\tnu_0<1$, $c=-1/2$, $\tlambda_n=1/2$  \\[0.2cm]
  $q<-1$   & BFA  & $\tnu_0\in Z$, $c=q^{-\tnu_0-1} [\tnu_0+1]_q$, $\tlambda_n =
                       [n-1]_q$  \\[0.2cm]
               & UB   & $0\le\tnu_0<1$, $c=(q-1)^{-1}$, $\tlambda_n = (1-q)^{-1}$
\\[0.2cm]
  \hline 
\end{tabular}
\end{table}
%
%
\normalsize
\subsection{The Chaturvedi-Srinivasan oscillator algebra}
It is characterized by $q=1$ on the left-hand side of (\ref{2}), and 
$G(N)=q^N$ on the right-hand side, i.e.,$ [ a, \s ] = q^N $ 
\cite{chaturvedi}. $ F(N) $ is therefore solution of 
\be
F(x+1) -  F(x) = G(x) \: , \; \;  x \in R \; \; 
                                {\rm and} \; \; F(0) = 0 \, , 
\ee
and
\be
C= F(N) - \s a \, . 
\ee
We obtain 
\be 
F(N) = [N]_q \, , 
\ee
again, but now
\be
C = [N]_q - \s a \, , 
\ee
so that 
\be
\la_n = \frac{q^{\nu_0 + n}-1}{q-1} - c \, .
\ee 
The unirreps of this algebra are given in table 2. The classically 
singular representations are now the UB unirreps 
with $ q > 1 $, contrary to what happens for the Arik-Coon oscillator.
%
%
\footnotesize
\begin{table}[htb]
\caption{Unirrep classification for the Chaturvedi-Srinivasan oscillator algebra.}
\vspace{.5cm}
\begin{tabular}{lll}          
  \hline\\[-0.2cm] 
  $q$ & Type & Characterization \rule[-0.3cm]{0cm}{0.6cm}\\[0.2cm]
  \hline\\[-0.2cm] 
  $q>1$       & BFB & $\tnu_0\in R$, $c=[\tnu_0]_q$, $\tlambda_n = q^{\tnu_0} [n]_q$
                       \\[0.2cm]
                 & UB   & $0\le\tnu_0<1$, $c\le - (q-1)^{-1}$, $\tlambda_n = [\tnu_0+n]_q
                       - c$\\[0.2cm]
  $0<q<1$   & BFB & $\tnu_0\in R$, $c=[\tnu_0]_q$, $\tlambda_n = q^{\tnu_0}
                       [n]_q$\\[0.2cm]
  $-1<q<0$ & BFB & $\tnu_0\in2 Z$, $c=[\tnu_0]_q$, $\tlambda_n = q^{\tnu_0}
                       [n]_q$\\[0.2cm]
  $q=-1$    & FD  & $\tnu_0\in 2Z$, $p=1$, $c=0$, $\tlambda_n = \left(1 - 
                       (-1)^n\right)/2$ \\[0.2cm]
                & UB   & $\tnu_0=0$, $c<0$, $\tlambda_n = - c + \left(1 - (-1)^n\right)/2$
                       \\[0.2cm]
  $q<-1$   & BFA  & $\tnu_0\in2Z+1$, $c=[\tnu_0+1]_q$, $\tlambda_n = 
                       q^{\tnu_0+1} [n-1]_q$  \\[0.2cm]
  \hline 
\end{tabular}
\end{table}
%
%
\normalsize
\subsection{The Tamm-Dancoff oscillator algebra}
It is characterized by $ G(N) = q^N $ \cite{odaka}, so that
\be 
F(N) = q^{N-1} N \; \; \; \; \; \; {\rm and} \; \; \; \; \; \; 
   C = q^{-1} N - q^{-N} \s a \, , 
\ee
and
\be
\la_n = q^{\nu_0 + n - 1} (\nu_0 + n -q c ) \, .
\ee
Only one class of unirreps exists for this GDOA, as showed in table 3.
\footnotesize
%
%
\begin{table}[htb]
\caption{Unirrep classification for the Tamm-Dancoff oscillator algebra.}
\vspace{.5cm}
\begin{tabular}{lll}          
  \hline\\[-0.2cm] 
  $q$ & Type & Characterization \rule[-0.3cm]{0cm}{0.6cm}\\[0.2cm]
  \hline\\[-0.2cm] 
  $0<q\ne1$       & BFB & $\tnu_0\in R$, $c=q^{-1}\tnu_0$, $\tlambda_n = 
             q^{\tnu_0+n-1} n$ \\[0.2cm]
  \hline 
\end{tabular}
\end{table}
%
%
\normalsize 
\section{Conclusion}
We developed the representation theory of GDOA's. The classification of their
unirreps can be most easily performed in terms of the eigenvalues of the 
Casimir operators $U$ and $C$. We showed that the spectrum of the 
number operator 
$N$ is discrete. Under the assumption that it is nondegenerate, the unirreps
can fall into one out of four classes (BFB, BFA, FD, UB), bosonic, and 
fermionic or parafermionic Fock-space representations occuring as special 
cases.

We provided examples for each of these classes, althought in the FD case, 
only two-dimensional unirreps were encountered. Higher-dimensional 
FD unirreps do however arise for some known deformed oscillator 
algebras~\cite{cqb}.

Applications of deformed oscillator algebras have been restricted up to now to
their Fock-space representations. Whether non-Fock-space representations, such as
those considered here, may have some useful applications remains an
interesting open question.

{\it Note:} During the Colloquium, Mich\`ele Irac-Astaud pointed out to
us the existence of some previous related works on the classification
of GDOA's unirreps \cite{irac}, confirming the
present results.


\begin{thebibliography}{99}
%
\bibitem{arik} Arik M. and Coon D. D.: J. Math. Phys. {\it 17} (1976) 
524.
%
\bibitem{kuryshkin} Kuryshkin V.: Ann. Fond. L. de Broglie {\it 5} 
(1980) 111;\\
                     Biedenharn L. C.: J. Phys. A {\it 22} (1989) L873;\\
                     Macfarlane A. J.: J. Phys. A {\it 22} (1989) 4581.
%
\bibitem{jannussis} Jannussis A., Brodimas G., and  Mignani R.: 
J. Phys. A {\it 24} (1991) L775;\\
                     Jannussis A.: J. Phys. A {\it 26} (1993) L233;\\
                     McDermott R. J. and  Solomon A. I.: J. Phys. A 
{\it 27} (1994) L15;\\
                     Meljanac S.,  Milekovi\'c M., and  Pallua S.:   
Phys. Lett. B {\it 328} (1994) 55.
%
\bibitem{daska91}  Daskaloyannis C.: J. Phys. A {\it 24} (1991) L789;\\
                     Bonatsos D. and  Daskaloyannis C.: Phys. Lett. B 
{\it 307} (1993) 100.
%
\bibitem{katriel}  Katriel J. and  Quesne C.: J. Math. Phys. {\it 37} 
(1996) 1650.
%
\bibitem{kulish}  Kulish P. P.: Theor. Math. Phys. {\it 86} (1991) 108;\\
                  Rideau G.: Lett. Math. Phys. {\it 24} (1992) 147;\\
                  Chaichian M.,  Grosse H., and  Pre\v snajder P.: 
J. Phys. A {\it 27} (1994) 2045.
%
\bibitem{gdoarep} Quesne C. and Vansteenkiste N.: Helv. Phys. Acta 
in press.
%
\bibitem{oh}  Oh C. H. and  Singh K.: J. Phys. A {\it 27} (1994) 5907.\\
              Quesne C. and  Vansteenkiste N.: J. Phys. A {\it 28} 
(1995) 7019.
%
\bibitem{jordan}  Jordan T. F.,  Mukunda N., and  Pepper S. V.: 
J. Math. Phys. {\it 4} (1963) 1089.
%
\bibitem{aizawa}  Aizawa N.: Phys. Lett. A {\it 177} (1993) 195;\\
                  Delbecq C. and  Quesne C.: Phys. Lett. B {\it 300} 
(1993) 227.
%
\bibitem{chaturvedi}  Chaturvedi S. and  Srinivasan V.: Phys. Rev. A 
{\it 44} (1991) 8020.
%
\bibitem{odaka}  Odaka K.,  Kishi T., and  Kamefuchi S.: J. Phys. A 
{\it 24} (1991) L591;\\
                 Chaturvedi S.,  Srinivasan V., and  Jagannathan R.: 
Mod. Phys. Lett. A {\it 8} (1993) 3727.
%
\bibitem{cqb}  Quesne C.: Phys. Lett. A {\it 193} (1994) 245.
%
\bibitem{irac} Irac-Astaud M. and Rideau G.: On the existence of
quantum bihamiltonian systems: the harmonic oscillator case, Preprint 
PAR-LPTM 92; Lett. Math. Phys. {\it 29} (1993) 197; Deformed quantum
harmonic oscillator, Proc. Third Int. Wigner Symposium (Oxford, 1993),
to appear; Theor. Math. Phys. {\it 99} (1994) 658.
%

\end{thebibliography}
\end{document}